\RequirePackage[leqno]{amsmath}
\documentclass[leqno]{amsart}
\usepackage{enumitem}
\setlist[itemize]{topsep=0pt,after=\vspace{1.5\baselineskip}}
 \usepackage[colorlinks=true]{hyperref}
\usepackage[T1]{fontenc}
\usepackage[margin=0.5in]{geometry}

\setlist[itemize]{noitemsep, topsep=0pt}

\def\R{\mathbb R} \def\N{\mathbb N} 
\def\id{\int_{\R^n}} 
\def\igrad{ \lvert \nabla \rho^\frac{m+p-1}{2}\rvert^2}
 
\def\R{\mathbb R} \def\N{\mathbb N} 
\def
\@cite
#1#2{[{{\bfseries #1}\if@tempswa , #2\fi}]}
\newtheorem{theorem}{Theorem}[section]
\newtheorem{corollary}[theorem]{Corollary}
\newtheorem{lemma}[theorem]{Lemma}

\newtheorem{remark}{Remark}

\title[Uniformly bounded solutions to an attraction-repulsion
model] %Use the shortened version of the full title
     {Influence of nonlinear productions on the global solvability of  an attraction-repulsion chemotaxis system}
\author[]{}
\title[Boundedness of solutions to a  nonlinear aggregation-diffusion equation] %Use the shortened version of the full title
      {%A nonlocal reaction chemotaxis model with arbitrarily positive diffusion exponent
Boundedness for a nonlocal reaction chemotaxis model even in the attraction-dominated regime}
\author[T. Li and G. Viglialoro]{}
\keywords{Aggregation-diffusion equations, Nonlinear parabolic equations, Global existence, Keller--Segel models, Nonlocal reactions.
\textit{$^*$Corresponding author}: giuseppe.viglialoro@unica.it}
\begin{document}
\maketitle
\maketitle

\centerline{\scshape Tongxing Li$^\natural$ \and Giuseppe Viglialoro$^{\sharp,*}$}
\medskip
{
% \footnotesize
% \centerline{$^1$Institut für Mathematik}
% \centerline{Universität  Paderborn }
% \centerline{Warburge
%r  Str.    100,  33098  Paderborn (Germany)}
  \centerline{$^\natural$School of Control Science and Engineering}
 \centerline{Shandong University}
 \centerline{Jinan, Shandong, 250061. (P. R. China)}
 \medskip
 \centerline{$^{\sharp}$Dipartimento di Matematica e Informatica}
 \centerline{Universit\`{a} di Cagliari}
 \centerline{Via Ospedale 72, 09124. Cagliari (Italy)}
%  \centerline{E-mail: giuseppe.viglialoro@unica.it}
 \medskip
}
\bigskip
\begin{abstract}
This work deals with a parabolic chemotaxis model with nonlinear diffusion and nonlocal reaction source. The problem is formulated on the whole space and, depending on a specific interplay between the coefficients associated to  such diffusion and reaction, we establish that all given solutions are uniformly bounded in time. To be precise, we study these
attractive (sign ``$+$'') and repulsive (sign ``$-$'') following models, formally described by the Cauchy problems
\begin{equation}\label{problem_abstract}
\tag{$\Diamond$}
\begin{cases}
\rho_t=\Delta \rho^m \pm \nabla \cdot \Big(\rho \nabla \Big(\frac{|x|^{2-n}}{2-n}*\rho\Big)\Big)+a\rho^\eta-b\rho^\alpha\id \rho^\beta dx   &  x\in \R^n, t\in (0,T_{\max}), \\
\rho(x,0)=\rho_{0}(x) & x\in  \R^n,
\end{cases}
\end{equation}
for $n\geq 3$, $m,a,b,\alpha,\eta>0$ and $\beta\geq 1$. By denoting with $T_{\max}$ the maximum time of existence of any nonnegative weak solution $\rho$ to problems \eqref{problem_abstract}, we prove that despite any large-mass initial data $\rho_0$, for any $\eta>0$ and arbitrarily small diffusive parameter $m>0$, whenever $\alpha+\beta$ surpasses some computable expression depending on $m, \eta$ and $n$, $T_{\max}=\infty$ and $\rho$ is uniformly bounded.

On the one hand, this paper is in line with claims established for $a=b=0$, where the same conclusion holds true in, respectively:
\begin{itemize}
	\item the repulsive scenario, under the assumption $m>0$ (adaptation of the case $m>1-\frac{2}{n}$, in Carrillo and Wang \cite{CarrilloWang});
	\item the attraction scenario, under the assumption $\frac{2n}{n+2}<m<2-\frac{2}{n}$ and for small initial data (Chen and Wang in \cite{ChenWangCriterionBlowup}).
\end{itemize}
On the other, for the attractive case with $a=b=m=1$ and $\eta=\alpha$,  this investigation also extends a result derived by Bian, Chen and Latos in \cite{BianChenEvangelosLinearNonlocal}.
\end{abstract}
\section{Introduction and motivations}\label{IntroSection}
This paper is, mostly, dedicated to the following Cauchy problem
\begin{equation}\label{problem}
\begin{cases}
\rho_t=\Delta \rho^m + \nabla \cdot (\rho \nabla (U*\rho))+a\rho^\eta-b\rho^\alpha\id \rho^\beta dx   &  x\in \R^n, t\in (0,T_{\max}), \\
%v_t=\Delta v-u^\beta v &  x\in \Omega, t>0, \\
%\frac{\partial u}{\partial \nu}=\frac{\partial v}{\partial \nu}=0 & x\in \partial \Omega, t>0, \\
\rho(x,0)=\rho_{0}(x) & x\in  \R^n,
\end{cases}
\end{equation}
for $n\geq 3$, $m>0$, $a,b,\alpha,\eta>0$, $\beta\geq 1$,  $\rho_0(x)\in L^1_+(\R^n)\cap L^\infty(\R^n)$, $U\in L^1_{\textrm{loc}}(\R^n)$ and $T_{\max}$ denoting the maximum time up to which its solutions are defined. The partial differential equation appearing above generalizes 
\begin{equation}\label{problemOriginalKS}
\rho_t=\Delta \rho + \nabla \cdot (\rho \nabla (U*\rho))\quad   x\in \R^n, t\in (0,T_{\max}),
\end{equation}
proposed in the pioneer papers by Keller and Segel (\cite{K-S-1970,Keller-1971-MC,K-S-1971}) to model the dynamics of populations (as for instance cells or bacteria) 
arising in mathematical biology. Precisely, by indicating with  $\rho = \rho(x, t)$ a certain cell density at the position $x$ and at the time $t$, the equation describes how a given chemotactical impact of the chemosensitivity (power-law potentials) $U(x)$, possibly justified by the presence of a chemical signal, may break the natural diffusion  (associated to the Laplacian operator, $\Delta \rho$) of the cells, initially organized accordingly to a certain configuration (the law $\rho_0(x)=\rho(x,0)$ in \eqref{problem}),  and even strongly influence their motion, leading the system to its chemotactic collapse (blow-up at finite time with appearance of $\delta$-formations). In the literature there are many contributions dedicated to the comprehension of this phenomena, especially in two-dimensional settings, and when the so-referred attractive Newtonian potential $U(x)=+\log(|x|)$ is fixed. In this regard, in  \cite{BlanchetCalvezCarrillo,BlanchetCarrilloNader,DOLBEAULT2004,JaLu,Nagai} the interested reader can find an extensive and rigorous theory on existence and properties of global, uniformly bounded or blow-up solutions to the initial problem associated to \eqref{problemOriginalKS}, especially in terms of the initial mass of the cell distribution, i.e.,  $M=\id \rho_0dx.$ Indeed, the mass of the bacteria, preserved in time for this model, appears as a critical parameter (see, for instance,  \cite{BlanchetDolbeaultPerthame2006} and \cite{ChildressPercus1981}); more exactly the value $M_c=8\pi$ establishes that when the diffusion overcomes the self-attraction ($M<M_c$), global in time solutions are expected, whereas when the self-attraction dominates the diffusion ($M>M_c$), blow-up solutions at finite time may be detected. (As to the corresponding asymptotic analysis we refer to \cite{BlancheCarlenCarrillo2012_JFA,DiFrancescoDolak-StrussSIAM2006,CamposDolbeault2014}.)

Besides the largeness of their initial distribution, the chemotactic behavior the cells toward their self-organization, may be sensitively influenced by other factors; the impacts of the diffusion (weaker or stronger) and of the chemosensitivity (attractive, repulsive or both), as well as the presence of external sources affecting the cells' density. Model \eqref{problem} is an example that combines all these aspects, exactly as specified. The parameter $m>0$ highlights the nonlinearity effect of the diffusion of the population; essentially, the larger $m$ the higher is the repulsion between each cell, inducing the system to a natural equilibrium. On the other hand, the presence of the first local addendum in the reaction term (corresponding to an increase of the population) and the second nonlocal addendum (corresponding to a decrease of the population) indicates that cells are, respectively, produced and/or consumed throughout time. (We suggest, primarily, \cite[Introduction 0.]{Souplet_Gradient}  and, also,  \cite[Remark 1.]{LiPintusViglialoroZAMP} for a discussion on the role of external sources in reaction-diffusion problems.) As to the description of the taxis-driven behavior of the cells, herein we are mainly interested in models describing cell aggregations and coalescence. For this reason we consider only fully attractive potentials (this choice is justified by technical reasons, as explained in item  \ref{NonNewtonianOpenProblems} of $\S$\ref{OpenProblemsSection}. below), expressed by
\begin{equation}\label{ExpressionPotentialU}
U(x)=\frac{|x|^{2-n}}{2-n}, \quad n\geq 3,
\end{equation}
a particular case of  
\begin{equation}\label{ExpressionPotentialUGeneral}
V(x)=\lambda_1\frac{|x|^A}{A}-\lambda_2 \frac{|x|^B}{B},  \quad \lambda_1,  \lambda_2 \in \{0,1\} \textrm{ and proper } A \textrm{ and }  B.
\end{equation}
Unlike $U$, $V$ is especially appropriate to idealize interaction of the cells' distribution influenced by both attractive and repulsive agents: attractive-repulsive for $\lambda_1=\lambda_2=1$, fully attractive for $\lambda_1=1$ and $\lambda_2=0$ and fully repulsive for $\lambda_2=1$ and $\lambda_1=0$. The attractive Newtonian potential \eqref{ExpressionPotentialU} employed in this investigation, naturally comes from \eqref{ExpressionPotentialUGeneral} for the triplet $(\lambda_1,\lambda_2,A)=(1,0,2-n)$. Oppositely, $(\lambda_1,\lambda_2,A)=(0,1,2-n)$ characterizes a fully repulsive taxis-driven term, and by replacing $U(x)$ in \eqref{problem} with $-U(x)$ it is rather conceivable that the resulting model prevents gathering phenomena and provides more smoothness and equilibrium to the system more efficiently than the case with fully attractive potential. 

As to some results concerning the blow-up suppression for cells densities involved in chemotaxis model, in the recent paper \cite{CarrilloWang}, this uniform-in-time $L^\infty$ criterion for weak solutions $\rho$ to system \eqref{problem}, with $n\geq2$, $a=b=0$ and $U(x)=V(x)$ as in \eqref{ExpressionPotentialUGeneral} and any $\rho_0\in L^1_+(\R^n)\cap L^\infty(\R^n)$, is derived: 
$$
\textrm{For \;}
\begin{cases}
(\lambda_1,\lambda_2)=(1,1),\, 2-n\leq B<A\leq 2  \textrm{ and }  m>1-\frac{2}{n}\\
(\lambda_1,\lambda_2)=(1,0),\, 2-n<A\leq 2 \textrm{ and }  m>1-\frac{A}{n}
\end{cases}
\;\; \rho \in L^\infty(\R_+;L^\infty(\R^n)).
$$
Indeed, the fully attractive Newtonian potential is deeply discussed in  \cite{ChenWangCriterionBlowup} and \cite{LiuWang}. In these researches the authors established, inter alia, this exact criterion for global existence and blow-up of solutions to \eqref{problem}, in the case $n\geq 3$, $a=b=0$,  $m\in (\frac{2n}{n+2},2-\frac{2}{n})$ and $U(x)$ as in \eqref{ExpressionPotentialU}: There exists $s^*=s^*(m,n) >0$ such that for any $\rho_0\in L^1_+(\R^n)\cap L^m(\R^n)$ 
$$
\begin{cases}
\textrm{if \;} \lVert \rho_0 \rVert_{L^\frac{2n}{n+2}(\R^n)}<s^*,  \textrm{ then problem \eqref{problem} has a globally bounded solution,}\\
\textrm{if \;}\lVert \rho_0 \rVert_{L^\frac{2n}{n+2}(\R^n)}>s^*,  \textrm{ then problem \eqref{problem} has a blow-up solution.}
\end{cases}
$$
Unlike the linear diffusion case \eqref{problemOriginalKS}, where the critical mass $M_c$ is $n$-independent, the above criterion implies that the size of the initial mass may have no crucial role on the existence of global or local in time solutions to nonlinear diffusion chemotaxis-systems. The key factor is given by some specific interplay between the diffusion coefficient $m$ and the dimension $n$; this is especially observed at high dimensions, where larger values of the diffusion parameters are required to ensure globability and boundedness of solutions. (This phenomenon holds true also in Keller--Segel systems defined in bounded domains; see \cite{HashiraIshidaYokotaJDEBlow-UP},  \cite{MarrasNishinoViglialoro} and \cite{WinDj}.)
\section{Presentation of the main results and open problems}
\subsection{Claims of the main results} 
In accordance to what discussed above, we are interested in the introduction of the external source $a\rho^\eta-b\rho^\alpha\id \rho^\beta dx$ exactly in the perspective toward blow-up prevention  in model \eqref{problem} without invoking any largeness assumption on $m$. Since as mentioned the difference in the source  stands for a competition between a birth contribution, favoring instabilities of the species (especially for large values of $a$ and/or $\eta$), and a death one opportunely contrasting this instability (especially for large values of $b$, $\alpha$ and $\beta$), some questions seem natural. For the  attractive Newtonian potential case, can one expect that the external source suffices to enforce globability of solutions, even for any arbitrarily small and $n$-independent positive value the diffusion parameter $m$? Are, conversely, some restrictions on $m$, $n,$ and/or $a, b$, $\alpha, \beta$ and $\eta$ required? And when the repulsive Newtonian potential is considered, may these restrictions (as expected) be weakened? Our main results positively address these questions in the sense that \textit{despite any fixed small value of the diffusion (i.e., $m$) and dampening  (i.e., $b$) parameters, any largeness of the initial data (i.e., $\rho_0$) and growth coefficient and power (i.e., $a$ and $\eta$), by sufficiently magnifying the impact of the powers associated to the death rate  (i.e., $\alpha$ and $\beta$)  of the source provides uniform-in-time boundedness of solutions to model \eqref{problem}.}

What said is formally claimed as follows.

% , let us first fix these mutual relations on the parameters $\alpha,\beta,\eta,m>0$,  also in term of $l:=\frac{2n}{n-2}$, a precise value related to the space dimension $n$, which herein we suppose to satisfy $n\geq 3$: 
%\begin{align}
%\label{b}\tag{\textsf{\textsf{H1}}}
%\alpha+\beta>\max\{(4 l - 4 - m l)/(p - 2), (2 \eta l - m l - 2 \eta)/(l - 
%    2) , m , 2 - m l/(l - 2)\}.
%\end{align}
%
%\begin{align}
%\label{b1}\tag{\textsf{\textsf{H2}}}
%\alpha+\beta>\max\{(4 l - 4 - m l)/(p - 2), (2 \eta l - m l - 2 \eta)/(l - 
%    2) , 2 - m l/(l - 2)\}.
%\end{align}
%These represent exactly our main assertions concerning the attractive-repulsive case ($\lambda=1$): 
\begin{theorem}[Attractive Newtonian potential]\label{MainTheorem}   
%{\textcolor{blue}{
For any $x\in \R^n$, with $n\geq 3$, let $U(x):=\frac{|x|^{2-n}}{2-n}$ and $\rho_0 \in L^1_+(\R^n)\cap L^\infty(\R^n)$. Moreover, let $l:=\frac{2n}{n-2}$, $a,b>0$, $\beta\geq 1$, $m,\alpha$ and $\eta$ be positive reals complying with \begin{align}
\label{b}\tag{\textsf{\textsf{H1}}}
\alpha+\beta>\max\Big\{\frac{4 l - 4 - m l}{l-2},\frac{2 \eta l - m l - 2 \eta}{l - 
    2}, m\Big\}.
\end{align}
If $\rho$ is a nonnegative and local weak solution to problem \eqref{problem} on $\R^n\times (0,T_{\max})$, then $T_{\max}=\infty$ and there exists a $C>0$ such that 
\begin{equation*}
\lVert \rho (\cdot,t)\rVert_{L^\infty(\R^n)}\leq C\quad \textrm{for all }t >0.
\end{equation*}
%\begin{equation}\label{GlobalBoundednessuinftyvW1infty}
%\lVert u(\cdot,t) \rVert_{L^\infty(\Omega)}\leq C \quad \textrm{for all} \quad t\in(0,\infty).
%\end{equation}
\end{theorem}
When the taxis-driven term is replaced by the  repulsive Newtonian potential, model \eqref{problem} is clearly more inclined to a natural smoothness under less restrictive assumptions than the attraction case. 
\begin{corollary}[Repulsive Newtonian  potential]\label{MainTheoremCorollary}   
%{\textcolor{blue}{
For any $x\in \R^n$, with $n\geq 3$, let $U(x):=-\frac{|x|^{2-n}}{2-n}$ and $\rho_0 \in L^1_+(\R^n)\cap L^\infty(\R^n)$. Moreover, let $l:=\frac{2n}{n-2}$, $a,b>0$, $\beta\geq 1$, $m,\alpha$ and $\eta$ be positive reals complying with 
\begin{align}
\label{b1}\tag{\textsf{\textsf{H2}}}
\alpha+\beta>\max\Big\{\frac{2 \eta l - m l - 2 \eta}{l - 
    2},m\Big\}.
\end{align}
 If $\rho$ is a nonnegative and local weak solution to problem \eqref{problem} on $\R^n\times (0,T_{\max})$, then $T_{\max}=\infty$ and there exists a $C>0$ such that 
\begin{equation*}
\lVert \rho (\cdot,t)\rVert_{L^\infty(\R^n)}\leq C\quad \textrm{for all }t >0.
\end{equation*}
%\begin{equation}\label{GlobalBoundednessuinftyvW1infty}
%\lVert u(\cdot,t) \rVert_{L^\infty(\Omega)}\leq C \quad \textrm{for all} \quad t\in(0,\infty).
%\end{equation}
\end{corollary}
\begin{remark}
In accordance with the nomenclature and the analysis from the interesting paper \cite{CarrilloHoffmannMaininiVolzone2018} (see also \cite{CalvezCarrilloHoffmannFairCompetition}, \cite{CalvezCarrilloHoffmannFairCompetitionTheGeometry} and correlated references therein), by introducing the fair-competition
regime associated to the attractive Newtonian potential, given by $m_c=1-\frac{2-n}{n}$, Theorem \ref{MainTheorem} establishes that if $\alpha+\beta$ is larger than some value depending on $m, n$ and $\eta$, then all solutions are global also for $m>0$ infinitesimally close to 0, and hence even in  the so-called  attraction-dominated regime, corresponding to the range $0 < m < m_c$. 
\end{remark}
\begin{remark}
As it will be discussed in the lines of the paper, the uniform-in-time estimate of the cells' density $\rho(x,t)$ in the $L^\infty$-norm passes through that in the $L^p$-norm, for some $p>1$. (a) For $a=b=0$ and $U(x)=-\frac{|x|^{2-n}}{2-n}$ in problem \eqref{problem}, such an estimate is immediately obtained; indeed, if in \cite[(2.32) of Theorem 2.2.]{CarrilloWang}  we refrain by considering the contribution from the attractive potential, the assumption $m>1-\frac{2}{n}$ therein used is unnecessary, and $m>0$ is sufficient to provide estimate (2.31). On the contrary, the introduction of the external source, in particular of the addendum tied to the birth rate, makes the analysis technically more involved and, accordingly to Corollary \ref{MainTheoremCorollary} (as well as to biological intuitions), it is required some strong effect of the dampening impact to have global existence. (b) For $a=b=0$ and $U(x)=\frac{|x|^{2-n}}{2-n}$ in problem \eqref{problem}, it is seen from \cite[(1) of Theorem 1.1.]{ChenWangCriterionBlowup} that in order to guarantee boundedness of solutions in the absence of the external source, especially of the non-local part contrasting the growth of the cell population, small-mass initial data and a restricted range of possible values of the diffusion parameter are required. Consequently Theorem \ref{MainTheorem} is an extension of the aforementioned. 

In addition, for $a=b=m=1$ and $\eta=\alpha$, assumptions \eqref{b} are simplified into $\alpha<1+\frac{2\beta}{n}$, so recovering for $\sigma=1$ the corresponding first condition in \cite[Theorem 1.2.]{BianChenEvangelosLinearNonlocal}, which henceforth is a particular case of Theorem \ref{MainTheorem}. 
\end{remark}
\subsection{Some open problems: hints and difficulties}\label{OpenProblemsSection}
In the next items we discuss some open questions, giving for each of them some considerations.
\begin{enumerate}[label=\roman*)]
\item  \textit{Local-existence.} As to the solvability of  \eqref{problem}, herein we follow the same approach than in \cite{CarrilloWang}, i.e., from now on we will assume that any $\rho_0 \in L_+^1(\R^n)\cap L^\infty(\R^n)$ emanates a \textit{sufficiently smooth local weak solution} $\rho$
and we dedicate to derive some \textit{a priori} uniform-in-time $L^p$ estimates for $\rho$. (In particular we might tacitly avoid to explicitly mention $\rho_0$.) In particular, our investigation is focused on the question concerning the maximum time $T_{\max}$ of existence of such solutions, for which the following \textit{extensibility criterion} holds true (see, for instance, \cite[Theorem 1.1.]{SugiyamaDIE2006}):  
\begin{equation}\label{ExtensibilityCrit}
\begin{aligned}
& \triangleright)\; T_{\max}=\infty, \textrm{so that $\rho$ remains bounded for all $x \in\R^n$ and all time $t>0$}, \\
& \triangleright)\;  \textrm{$T_{\max}$ finite (the blow-up time),  so that  $\limsup_{t \rightarrow T_{\max}}\lVert \rho(\cdot,t)\rVert_{L^\infty(\R^n)} =\infty$.} 
\end{aligned}
\end{equation} 
We indicate to the reader \cite{BianLiu2013Existence,SugiyamaDIE2006}, where one can find the analysis concerning the  local-in-time existence  of weak solutions to problem \eqref{problem}  for $a=b=0$ and $m>2-\frac{2}{n}$, as well as \cite{BianChenEvangelosLinearNonlocal} where for $a=b=m=1$ and $\eta=\alpha$ the same question is addressed for strong solutions. Hereafter, a well-posedness theory might be developed by suitable adaptations of the mentioned contributions.  
\item  \textit{Pure logistic sources.}
The well-known \textit{logistic model} in population dynamics is described by the equation
\begin{equation}
\frac{d p}{d t}=ap-bp^2, \nonumber
\end{equation}
where $a,b>0$ are related to the carrying capacity of the species $p=p(t)$, also associated to its growth and death rate. Conversely to the simpler case
$$\rho_t=\Delta \rho^m + \nabla \cdot (\rho \nabla (U*\rho)) \quad t\in (0,T_{\max}), \quad  \, \rho(x,0)=\rho_0(x)$$
for which weak solutions are such that the conservation mass $M_0:=\id \rho_0(x)dx=\id \rho dx$ holds (see \cite{BianLiu2013Existence} for more details), adding the logistic term to the model produces 
$$\rho_t=\Delta \rho^m + \nabla \cdot (\rho \nabla (U*\rho))+a \rho -b \rho^2\quad t\in (0,T_{\max}), \quad  \, \rho(x,0)=\rho_0(x)$$
and a formal integration over $\R^n$ gives
\begin{equation*}
\frac{d}{d t}\id \rho (x,t) dx =a\id\rho(x,t) dx-b\id\rho^2(x,t)dx\quad t \in (0,T_{\max}).
\end{equation*}
If the above logistic model was defined in a bounded and smooth domain $\Omega$ of $\R^n$ (see, for instance, \cite{MarrasViglialoroMathNach} and \cite{ViglialoroWoolley2018}), the Young inequality would provide this bound for $\lVert \rho \rVert_{L^1(\Omega)}$,
\begin{equation*}
\id \rho dx \leq \max\Big\{\int_\Omega \rho_0dx,\frac{a|\Omega|}{b}\Big\}\quad \textrm{on } (0,T_{\max}),
\end{equation*}
a cornerstone toward uniform-in-time $L^p$-bounds of $\rho$ and, in turn, also in $L^\infty$. The lack of an upper bound for $\id \rho dx$ leads to some technical restrictions, so that the general machinery fails. Indeed, the nonlocal term in $a\rho^\eta-b\rho^\alpha\id \rho^\beta dx$ helps to circumvent the problem, since the time derivative of $\id \rho^p dx$, for proper $p>1$, produces an addendum of the type $\id \rho^{p+\alpha-1}dx \id\rho^\beta dx$. As specified in \cite{BianChenFujitaExponent} and \cite{BianChenEvangelosLinearNonlocal}, this is crucial in the treatment of certain contributions involving $\id \rho^p dx$. 
\item \label{NonNewtonianOpenProblems} \textit{Non-Newtonian potentials.} With the fundamental solution of the Laplace equation in our hands, it will be seen that $\Delta (U* \rho)$ is proportional to $\rho$. Consequently, considering the Newtonian powers allows that standard testing procedures and usages of functional inequalities smoothly work.  This proportionality does not apply for non-Newtonian potential, and the analysis resulting  from general expressions as in  \eqref{ExpressionPotentialUGeneral} appears more challenging.
\end{enumerate}
\section{Fixing some parameters and functional inequalities. Organization of the paper}
This set of lemmas will be used in our logical steps to prove the main results. We first recall the fundamental solution of Laplace's equation (Lemma  \ref{LemmaLaplacianULimitCase}),  and then summarize some general functional inequalities (Lemmas \ref{InterpolationInequalityLemma} and \ref{InequalityG-NLemma}); we note that even though the validity of these lemmas concerns a larger class of functions, we claim them directly for all weak solutions $\rho$ defined $\R^n \times (0,T_{\max})$ to problem \eqref{problem}. Moreover, by adjusting the data $m,\alpha,\beta,\eta$ and $n$ of the same problem, in Lemma \ref{Lemmapbarra} we fix the value $\bar{p}$ of an important parameter, so that the employment of further inequalities leads, for any $p>\bar{p}$, to uniform-in-time $L^p$-bound on $(0,T_{\max})$ of these solutions $\rho$. Thereafter, relying on the gained $L^p$-bound, an adaptation (Lemmas \ref{RecursiveInequalityLemma} and  \ref{Sequel_p_r_theta_Lemma}) of the bootstrap Moser--Alikakos iterative method (see \cite{AlikakosTechnique}) is used to prove uniform-in-time boundedness on $(0,T_{\max})$ ($\S$\ref{SectionProofsClaims}.). Finally, the extensibility criterion \eqref{ExtensibilityCrit} shows that the aforementioned local weak solutions are actually globally bounded weak solutions ($\S$\ref{AprioriEstimatesSection}.).
\begin{lemma}\label{LemmaLaplacianULimitCase}
For any $n\geq 3$, the potential $U(x)$ defined in \eqref{ExpressionPotentialU} is such that on $\R^n$  
\begin{equation}\label{LaplacianoUintheLimitcase}
\Delta U (x)=n\alpha_n\delta(x),
\end{equation}
$\delta(x)$ denoting the Dirac measure on $\R^n$ at the point 0 and where $\alpha_n=\frac{\pi^\frac{n}{2}}{\Gamma(\frac{n}{2}+1)}$ is the volume of the $n$-dimensional unit ball. 

In particular, for any local weak solution $\rho$  to problem \eqref{problem} defined in $\R^n\times (0,T_{\max})$ it holds that
\begin{equation}\label{LaplacianoConvolutioRhoUintheLimitcase}
\Delta (U* \rho) (x,t)=n\alpha_n\rho(x,t)\quad \textrm{on } \R^n\times (0,T_{\max}).
\end{equation}
\begin{proof}
It is known (see \cite[page 22]{Evans-2010-PDEs}) that for any $n\geq 3$ the function $\phi(x)=\frac{1}{n(2-n)\alpha_n} | x | ^{2-n}$ solves the Laplace equation $-\Delta \phi=\delta$, so  \eqref{LaplacianoUintheLimitcase} is obtained. In turn relation \eqref{LaplacianoConvolutioRhoUintheLimitcase} follows from convolution operations and the sifting property of the $\delta$-distribution. 
\end{proof}
\end{lemma}
\begin{lemma}\label{Lemmapbarra} 
	For $n\in \N$, with $n\geq 3$, let $l:=\frac{2n}{n-2}$, $\alpha,\eta,m>0$ and $\beta\geq 1$ satisfy the assumptions in \eqref{b}. Then by defining
		\begin{equation}\label{ConstantForTechincalInequality_Barp}
	\bar{p}:=\max 
	\begin{Bmatrix}
	1 \vspace{0.1cm}   \\ 
	(l+2-lm)/(l-2) \vspace{0.1cm}  \\
	-\alpha-\beta+	(5l-2lm-2)/(l-2) \vspace{0.1cm}\\
		(2 \eta+l-l m-2)/(l-2)  \vspace{0.1cm} \\
	-\alpha-\beta+	(l+2\eta l-2-2m l)/(l-2)  \vspace{0.1cm}\\
		m-1  \vspace{0.1cm}\\
		\alpha +\beta -1 \vspace{0.1cm}\\
	\alpha +\beta +1-2\eta\vspace{0.1cm}\\
	1-\alpha-\beta+	2l(\eta-m)/(l-2)\vspace{0.1cm}\\
		1-\alpha-\beta+	2l(2-m)/(l-2)\vspace{0.1cm}\\
		(2\eta+l-lm-2)/(l-2) 
%		\alpha +\beta +1-2\eta  \vspace{0.1cm} 
% 
%	\\ 2l(2-m)/(l-2)+1-\alpha-\beta \vspace{0.1cm} \\
%	(2-(m-1))/(l-2) \vspace{0.1cm} \\
%	\beta+1-\alpha \vspace{0.1cm}\\ 
%	m-2\eta+1 \vspace{0.1cm}  \\ 
%	
%	2l(\eta-m)/(l-2)+1-\alpha -\beta \vspace{0.1cm}\\
%	\alpha+\beta+1-2\eta \vspace{0.1cm} \alpha+\beta-3 
	\end{Bmatrix}
	,
	\end{equation}
	   for all $p>\bar{p}$ it holds that $p'=\frac{p+\alpha-1+\beta}{2}\in (\beta,p+\alpha-1)$ and that 
	\begin{equation}\label{Lamda_0_Lambda_1_Lambda_eta} 
	\begin{cases}
	\tilde{\Lambda}_0:=\displaystyle \frac{p (l-2)+l(m-1)}{(l-2)(\alpha+\beta+p)-3l+\alpha l+2ml} &\in (0,1)\\
	\tilde{\Lambda}_1:=\displaystyle\frac{ p (l-2)+(m-1) l-2}{(l-2)(\alpha+\beta+p)-5 l+2+2 m l} &\in (0,1)\\
	\tilde{\Lambda}_\eta:=\displaystyle\frac{2-2 \eta+p (l-2)+l(m-1)}{(l-2)(\alpha+\beta+p)+2+2ml-l-2\eta l} &\in (0,1).
	\end{cases}
	%\end{minipage}
	\end{equation}

%	there exists $\bar{p}$ such that for all $p > \bar{p}$ and  $p'=\frac{p+\alpha-1+\beta}{2}$ these relations hold 
%	\vskip0.1cm
%	\begin{subequations}
%		\begin{minipage}{0.3\textwidth}
%			\begin{align}
%			\label{a3}
%			p>\alpha+\beta-3
%			\end{align}
%		\end{minipage}
%		\begin{minipage}{0.3\textwidth}
%			\begin{align}
%			\label{a1}
%			p>
%			\end{align}
%		\end{minipage}
%		\begin{minipage}{0.4\textwidth}
%			\begin{align}
%			\label{a4}
%			p>\frac{2l}{l-2}(2-m)-\alpha-\beta+1
%			\end{align}
%		\end{minipage}
%		\vskip0.1cm
%		\begin{minipage}{0.3\textwidth}
%			\begin{align}
%			\label{a8bis}
%			p>\alpha+\beta+1-2\eta
%			\end{align}
%		\end{minipage}
%		\begin{minipage}{0.3\textwidth}
%			\begin{align}
%			\label{a9}
%			p>
%			\end{align}
%		\end{minipage}
%		\begin{minipage}{0.4\textwidth}
%			\begin{align}
%			\label{a9bis}
%			p>
%			\end{align}
%		\end{minipage}
%		\vskip0.1cm
%		\begin{minipage}{0.3\textwidth}
%			\begin{align}
%			\label{a20}
%			p>\alpha+\beta-1 
%			\end{align}
%		\end{minipage}
%		\begin{minipage}{0.3\textwidth}
%			\begin{align}
%			\label{a16}
%			p>\frac{l-lm}{l-2}
%			\end{align}
%		\end{minipage}
%		\begin{minipage}{0.4\textwidth}
%			\begin{align}
%			\label{RangeOfpPrime}
%			p'\in (\beta,p+\alpha-1).
%			\end{align}
%		\end{minipage}
%	\end{subequations}

	\begin{proof}
		Simple algebraic considerations give all the claims. 
	\end{proof}
\end{lemma}
\begin{lemma}[\protect{Interpolation inequality}]\label{InterpolationInequalityLemma}
	Under the hypotheses of Lemma \ref{Lemmapbarra}, let $\bar{p}$ and $p'$ be therein defined. If $\rho$ is a nonnegative local weak solution to \eqref{problem} on $\R^n\times (0,T_{\max})$, then for every $p>\bar{p}$ it satisfies 
\begin{equation}\label{InterpolationInequa}
\lVert \rho \rVert_{L^{p'}(\R^n)}\leq \lVert \rho \rVert^{a_1}_{L^{\beta}(\R^n)}\lVert \rho \rVert^{1-a_1}_{L^{p+\alpha-1}(\R^n)}, 
\quad \text{ on }  (0,T_{\max})\quad \textrm{ with } a_1=\frac{\frac{1}{p'}-\frac{1}{p+\alpha-1}}{\frac{1}{\beta}-\frac{1}{p+\alpha-1}}\in [0,1].
\end{equation}
\begin{proof}
See \cite[page 93]{BrezisBook}.
\end{proof}
\end{lemma}
\begin{lemma}[\protect{Gagliardo--Nirenberg-type inequality}]\label{InequalityG-NLemma} 
For $n \in \N$, with $n\geq 3$, let $l:=\frac{2n}{n-2}$, $1\leq r<q<l$ and $\frac{q}{r}<\frac{2}{r}+1+\frac{2}{l}$. Moreover let us set
\begin{equation}\label{ValueOfGammaInGagliardoNiren}
\lambda:=\frac{\frac{1}{r}-\frac{1}{q}}{\frac{1}{r}-\frac{1}{l}}\in (0,1), \;\,\;\gamma:=\frac{2(1-\lambda) q}{2-\lambda q}.
\end{equation}
Then,  for any $\varepsilon>0$ there exist $C_1(n)>0$ and $C(\varepsilon)=C_1(n)\varepsilon^\frac{-\lambda q}{2-\lambda q}$ such that, if $\rho$ is a nonnegative local weak solution to \eqref{problem} on $\R^n\times (0,T_{\max})$ it holds that
\begin{equation}\label{InequalityTipoG-N} 
\| \rho \|_{L^{q}(\R^n)}^q \leq \varepsilon \| \nabla \rho \|_{L^{2}(\R^n)}^2+C(\varepsilon)\|\rho \|_{L^{r}(\R^n)}^{\gamma}\quad \text{ on }  (0,T_{\max}).
\end{equation}
\begin{proof}
The proof is consequence of Sobolev inequalities and detailed arguments can be found in \cite[Lemma 2.]{BianChenFujitaExponent}.
\end{proof}
\end{lemma}
\begin{lemma}\label{RecursiveInequalityLemma}
For some $L>0$ and all $k\in\N_0$, let $M_k \subset [1,\infty)$ such that
\begin{equation}\label{RecursiveInequality}
M_k\leq L M_{k-1}^{\theta_k}\quad \textrm{ for all }\;k\in \N,
\end{equation}
for some $\theta_k \subset (0,\infty)$ fulfilling $\theta_k\leq 2$ for all $k\in \N$. Then 
\begin{equation*}%\label{LastInequalityIterationOnMk}
M_k\leq L^k\cdot M_0^{2^k} \quad \textrm{ for all } \;k\in \N.
\end{equation*}
\begin{proof}
The proof is an adaptation of \cite[Lemma 4.3.]{WinklerExponentailDecay}. 
\end{proof}
\end{lemma}
\begin{lemma}\label{Sequel_p_r_theta_Lemma}  
Under the same assumptions of Lemma \ref{Lemmapbarra}, let $\bar{p}$ be therein defined, $k\in \N$ and $p_k:=2^k+\bar{p}$. Moreover, for any $m>0$ let us consider these positive sequences: 
\begin{equation*}
%\begin{cases}
q_{1,k}=\frac{2 (p_k + 1)}{m + p_k - 1},\quad
q_{\eta,k}=\frac{2 (p_k + \eta-1)}{m + p_k - 1},\quad
q_{0,k}=\frac{2 p_k }{m + p_k - 1}\quad \textrm{and}\quad r_k=\frac{2(p_{k-1}+1)}{m + p_k - 1}.
%\end{cases}
\end{equation*}
Then the sequences $\gamma_{1,k}, \gamma_{\eta,k}$ and $\gamma_{0,k}$ computed through \eqref{ValueOfGammaInGagliardoNiren} after previous calculations of $\lambda_{1,k}, \lambda_{\eta,k}$ and $\lambda_{0,k}$, are such that 
\begin{equation*}%\label{InequalityOn_mukForSomeD}
\mu_{1,k}=\frac{\gamma_{1,k}}{r_k}\leq 2,\quad \mu_{\eta,k}=\frac{\gamma_{\eta,k}}{r_k}\leq 2\quad \textrm{and}\quad \mu_{0,k}=\frac{\gamma_{0,k}}{r_k}\leq 2 \quad \textrm{ for all } \,k\in \N.
\end{equation*}
\begin{proof}
From their definitions, some tedious but standard reasoning show that $\mu_{1,k}, \mu_{\eta,k}$ and $\mu_{0,k}$ are increasing. On the other hand, since $p_k \nearrow +\infty$, we have that  $q_{1,k}$, $q_{\eta,k}$, $q_{0,k}$ as well as $r_k$ converge to 2, whereas $\lambda_{1,k}$, $\lambda_{\eta,k}$, $\lambda_{0,k}$ are infinitesimal.  We conclude by noting that also $\mu_{1,k}, \mu_{\eta,k}$ and $\mu_{0,k}$ have $2$ as limit.
\end{proof}
\end{lemma}
\section{Some a priori estimates: deriving uniform-in-time $L^p$-bounds for $\rho$}\label{AprioriEstimatesSection}
In this section, by establishing for $\id \rho^p dx$, with $p>\bar{p}$, an absorptive differential inequality for $t\in (0,T_{\max})$,  we obtain a uniform-in-time bound in $L^p(\R^n)$ for the cell density $\rho$. This will be achieved for any $m>0$.
\begin{lemma}\label{LEmmaBoundinlP}
	Under the hypotheses of Lemma \ref{Lemmapbarra} let $\bar{p}$ and $p'$ be therein defined. If $\rho$ is any local and nonnegative weak solution to problem \eqref{problem} on $\R^n\times (0,T_{\max})$, then there exists a $K=K(p)>0$ such that for every $p>\bar{p}$
\begin{equation}\label{uniformBoundForLpNormSecondCase}
\lVert \rho (\cdot, t)\rVert_{L^p(\R^n)}\leq K\quad \textrm{for all } t \in (0,T_{\max}).
\end{equation}
\begin{proof}
Up to rigorous limiting processes involving standard cut-off functions, for any $p>\bar{p}$, by using the equation in \eqref{problem} and the integration by parts formula (twice in the integral involving the potential), we can compute 
\begin{equation*}%\label{FirstStepMainDerivation}
\begin{split}
\frac{d}{d t} \id \rho^p dx=&p \id \rho^{p-1}\rho_t dx=-mp(p-1)\id \rho^{p+m-3}\lvert \nabla \rho\rvert^2 dx+(p-1)\id \rho^p \Delta (U*\rho)dx \\ &
\quad +p a\id \rho^{p+\eta-1}dx-pb  \id \rho^{p+\alpha-1}dx\id \rho^\beta dx\quad \textrm{on } (0,T_{\max}).
\end{split}
\end{equation*}
As to $\Delta (U*\rho)$, we invoke Lemma \ref{LemmaLaplacianULimitCase},  write $(p-1)\id \Delta (U*\rho) dx=n(p-1)\alpha_n\id \rho dx$ and take into consideration this pointwise identity 
$$\lvert \nabla \rho^\frac{m+p-1}{2}\rvert^2=2c_1 \rho^{m+p-3}|\nabla \rho|^2, \textrm{ with }c_1=\frac{2mp(p-1)}{(m+p-1)^2}\;\textrm{ for all } t\in (0,T_{\max}).$$ 
We hence arrive for all  $t\in(0,T_{\max})$ at 
\begin{equation}\label{FirstStepMainDerivation} 
\begin{split}
\frac{d}{d t} \id \rho^p dx&  +2c_1\id \lvert \nabla \rho^\frac{m+p-1}{2}\rvert^2dx+pb  \id \rho^{p+\alpha-1}dx\id \rho^\beta dx 
= p a\id \rho^{p+\eta-1}dx +n (p-1)\alpha_n \id \rho^{p+1}dx.
\end{split}
\end{equation}
Now we turn our attention to estimate the two contributions on the right-hand side and $\id \rho^p dx$ of the above relation in terms of the nonlocal on the left. This will be possible thanks to Lemmas \ref{InterpolationInequalityLemma} and \ref{InequalityG-NLemma}, and by collecting these bounds we will provide the desired absorptive differential inequality. 
\subsection*{\quad Estimating $n (p-1)\alpha_n \id \rho^{p+1}$.} Setting $q=2\frac{p+1}{m+p-1}$ and $r=\frac{2p'}{m+p-1}$, we have by relying on assumptions \eqref{b}, and in view of the properties of $\bar{p}$ defined in \eqref{ConstantForTechincalInequality_Barp} of Lemma \ref{Lemmapbarra},  that
$$1\leq r <q<l\quad \textrm{and}\quad\frac{q}{r}<\frac{2}{r}+1-\frac{2}{l}.$$ 
Subsequently, by using inequality \eqref{InequalityTipoG-N} with $\gamma$ (indicated now for reasons of clarity with $\gamma_1$) computed with relation \eqref{ValueOfGammaInGagliardoNiren}  precisely by using the above fixed values of $q$ and $r$, we can write on $(0,T_{\max})$ and for every $\varepsilon_1>0$ and some $C_1(\varepsilon_1)>0$
\begin{equation}\label{FirstEtimateIntp+1}
\begin{split}
n (p-1)\alpha_n  \id \rho^{p+1}dx&=n (p-1)\alpha_n \lVert \rho^\frac{q(m+p-1)}{2}\rVert_{L^q(\R^n)}^q\leq \varepsilon_1\id\igrad dx +C_1(\varepsilon_1)\bigg(\id \rho^{p'}dx\bigg)^\frac{\gamma_1(p+m-1)}{2p'}.
\end{split}
\end{equation} 
On the other hand, from the definition and the property of $p'$ in Lemma \ref{Lemmapbarra}, we can employ inequality  \eqref{InterpolationInequa}  so to have
\begin{equation}\label{EmploymentOfInterpolationonpPrime}
\bigg(\id \rho^{p'}dx\bigg)^\frac{1}{p'}=\lVert \rho \rVert_{L^{p'}(\R^n)}\leq \lVert \rho \rVert_{L^\beta(\R^n)}^{a_1}\lVert \rho \rVert_{L^{p+\alpha+1}(\R^n)}^{1-a_1}\quad \textrm{for all } t\in(0,T_{\max}).
\end{equation}  	
 Next, with a view to Lemma \ref{Lemmapbarra}, for $\Lambda_1=\frac{\gamma_1(p+m-1)}{2}$ and   $a=\frac{\frac{1}{\beta}-\frac{1}{p'}}{\frac{1}{\beta}-\frac{1}{p+\alpha-1}}$, it is seen, after some  tangled computations,  that  $\frac{a\Lambda_1}{p+\alpha-1}=\tilde{\Lambda}_1<1$ (recall \eqref{Lamda_0_Lambda_1_Lambda_eta}), and we get
\begin{equation*}
\lVert \rho \rVert_{L^{p'}(\R^n)}^{\Lambda_1}\leq (\lVert \rho \rVert_{L^\beta(\R^n)}^{\beta}\lVert \rho \rVert_{L^{p+\alpha-1}(\R^n)}^{p+\alpha-1})^{\frac{a\Lambda_1}{p+\alpha-1}}\lVert \rho \rVert_{L^{\beta}(\R^n)}^{\Lambda_1(1-a-\frac{a\beta}{p+\alpha-1})}\quad \textrm{on } (0,T_{\max}).
\end{equation*}
Finally, the Young inequality (note that $\lVert \rho \rVert_{L^{\beta}(\R^n)}^{\Lambda_1(1-a-\frac{a\beta}{p+\alpha-1})}=1$ due to $1-a-\frac{a\beta}{p+\alpha-1}=0$) gives for any $\delta_1>0$  a computable and positive $D_1(\delta_1)$ such that
\begin{equation}\label{FinalEstimateIntp+1}
C_1(\varepsilon_1)\bigg(\id \rho^{p'}dx\bigg)^\frac{\gamma_1(p+m-1)}{2p'}\leq \delta_1\id \rho^{p+\alpha-1}dx \id \rho^\beta dx+D_1(\delta_1)\quad \textrm{on } (0,T_{\max}).
\end{equation}
\subsection*{\quad Estimating $p a\id \rho^{p+\eta-1}dx$.} For $q=2\frac{\eta +p-1}{m+p-1}$
 and  $r=\frac{2p'}{m+p-1}$, again hypotheses \eqref{b} and  $\bar{p}$ in \eqref{ConstantForTechincalInequality_Barp} show that 
$$1\leq r <q<l\quad \textrm{and}\quad\frac{q}{r}<\frac{2}{r}+1-\frac{2}{l},$$ so having, through  \eqref{InequalityTipoG-N}, 
\begin{equation}\label{FirstEtimateIntp+eta+1}
\begin{split}
p a \id \rho^{\eta +p-1}dx&=p a  \lVert \rho^\frac{q(m+p-1)}{2}\rVert_{L^q(\R^n)}^q\leq \varepsilon_2\id\igrad dx+C_2(\varepsilon_2)\bigg(\id \rho^{p'}dx\bigg)^\frac{\gamma_\eta(p+m-1)}{2p'}\quad \textrm{for all } t\in (0,T_{\max}),
\end{split}
\end{equation}
for all  $\varepsilon_2>0$, some positive $C_2(\varepsilon_2)$ and $\gamma_\eta$  as in \eqref{ValueOfGammaInGagliardoNiren}.  Moreover, again for $a=\frac{\frac{1}{\beta}-\frac{1}{p'}}{\frac{1}{\beta}-\frac{1}{p+\alpha-1}}$ and thanks to \eqref{Lamda_0_Lambda_1_Lambda_eta}, if we consider $\Lambda_\eta=\frac{\gamma_\eta(p+m-1)}{2}$, we obtain that $\frac{a\Lambda_\eta}{p+\alpha-1}=\tilde{\Lambda}_\eta<1$. Subsequently, from \eqref{EmploymentOfInterpolationonpPrime}, we can deduce first
\begin{equation*}
\lVert \rho \rVert_{L^{p'}(\R^n)}^{\Lambda_\eta}\leq (\lVert \rho \rVert_{L^\beta(\R^n)}^{\beta}\lVert \rho \rVert_{L^{p+\alpha-1}(\R^n)}^{p+\alpha-1})^{\frac{a\Lambda_\eta}{p+\alpha-1}}\lVert \rho \rVert_{L^{\beta}(\R^n)}^{\Lambda_\eta(1-a-\frac{a\beta}{p+\alpha-1})}\quad \textrm{on } (0,T_{\max}),
\end{equation*}
and then, through Young's inequality,
\begin{equation}\label{FinalEstimateIntp+eta-1}
C_2(\varepsilon_2)\bigg(\id \rho^{p'}dx\bigg)^\frac{\gamma(p+m-1)}{2p'}\leq \delta_2\id \rho^{p+\alpha-1}dx \id \rho^\beta dx+D_2(\delta_2)\quad \textrm{on } (0,T_{\max}),
\end{equation}
for any $\delta_2>0$ and some computable $D_2(\delta_2)>0$. 
\subsection*{\quad Estimating $ \id \rho^{p}dx$.}
We use once more assumptions \eqref{b}, $\bar{p}$ and $a=\frac{\frac{1}{\beta}-\frac{1}{p'}}{\frac{1}{\beta}-\frac{1}{p+\alpha-1}}$. Then, for $q=\frac{2p}{m+p-1}$, $r=\frac{2p'}{m+p-1}$, we have 
$$1\leq r <q<l\quad \textrm{and}\quad\frac{q}{r}<\frac{2}{r}+1-\frac{2}{l},$$
and bound  \eqref{ValueOfGammaInGagliardoNiren} provides $\gamma_0$ such that due to relation \eqref{InequalityTipoG-N} for every $\varepsilon_3>0$ and some positive $C_3(\varepsilon_3)$ this inequality holds true:
\begin{equation}\label{FirstEtimateIntp}
\begin{split}
\id \rho^{p}dx&=\lVert \rho^\frac{q(m+p-1)}{2}\rVert_{L^q(\R^n)}^q\leq \varepsilon_3\id\igrad dx+C_3(\varepsilon_3)\bigg(\id \rho^{p'}dx\bigg)^\frac{\gamma_0(p+m-1)}{2p'}\quad \textrm{for all } t\in(0,T_{\max}).
\end{split}
\end{equation} 
On the other hand, $\frac{\gamma_0(p+m-1)}{2}=\Lambda_0$ and \eqref{Lamda_0_Lambda_1_Lambda_eta} entail $\frac{a \Lambda_0}{p+\alpha-1}= \tilde{\Lambda}_0<1$; consequently, 
relation \eqref{EmploymentOfInterpolationonpPrime} simplifies into 
\begin{equation*}
\lVert \rho \rVert_{L^{p'}(\R^n)}^{\Lambda_0}\leq (\lVert \rho \rVert_{L^\beta(\R^n)}^{\beta}\lVert \rho \rVert_{L^{p+\alpha-1}(\R^n)}^{p+\alpha-1})^{\frac{\Lambda_0}{p+\alpha-1}}\quad \textrm{on } (0,T_{\max}),
\end{equation*}
which for some couple $(\delta_3,D_3(\delta_3))$, both positive, reads through the Young inequality 
\begin{equation}\label{FinalEstimateIntp}
C_3(\varepsilon_3)\bigg(\id \rho^{p'}dx\bigg)^\frac{\gamma(p+m-1)}{2p'}\leq \delta_3\id \rho^{p+\alpha-1}dx \id \rho^\beta dx+D_3(\delta_3)\quad \textrm{on } (0,T_{\max}).
\end{equation}
\subsection*{\quad Providing the bound for $\id \rho^{p}dx$.}
The last step toward the proof of this lemma is considering for any $0<\varepsilon_3<2c_1$ in relation \eqref{FirstEtimateIntp},  $\varepsilon_1=\varepsilon_2=\frac{2c_1-\varepsilon_3}{2}$ in both \eqref{FirstEtimateIntp+1} and \eqref{FirstEtimateIntp+eta+1}, and $\delta_1=\delta_2=\delta_3=\frac{pb}{3}$ in \eqref{FinalEstimateIntp+1}, \eqref{FinalEstimateIntp+eta-1} and  \eqref{FinalEstimateIntp}. Plugging the results of these operations into  bound \eqref{FirstStepMainDerivation}, we obtain  that $ \id \rho^pdx $ has to satisfy for some positive $c_2$ (note that $\rho_0 \in L_+^1(\R^n)\cap L^\infty(\R^n)$ implies that $\rho_0 \in L^p(\R^n)$ for all $p\geq 1$)  this initial problem
\begin{equation*}\label{MainInitialProblemWithM}
\begin{cases}
\frac{d}{d t} \id \rho^p dx +  \id \rho^{p}dx
\leq c_2 \quad \textrm{for all } t\in(0,T_{\max}),\\
 \id \rho^p(x,0)dx= \id \rho^p_0(x)dx.
\end{cases}
\end{equation*}
An application of a comparison principle implies that
\begin{equation*}
\lVert \rho(\cdot,t)\rVert_{L^p(\R^n)}\leq \Big(\max\Big\{\id \rho^p_0(x)dx, c_2\Big\}\Big)^\frac{1}{p}:=K\quad \textrm{for all } t\in(0,T_{\max}).
\end{equation*}
%to obtain for $D=D_1(\frac{pb}{4})+D_2(\frac{pb}{4})$
%\begin{equation}\label{SecondStepMainDerivation}
%\begin{split}
%\frac{d}{d t} \id \rho^p dx&  +c_1\id \lvert \nabla \rho^\frac{m+p-1}{2}\rvert^2dx+\frac{pb}{2}  \id \rho^{p+\alpha-1}dx\id \rho^\beta dx 
%\leq D \quad \textrm{on } (0,T_{\max}).
%\end{split}
%\end{equation}
%Once bound \eqref{FirstStepMainDerivation}
%
%
%
%for $\varepsilon_3=c_1$ relation \eqref{FirstEtimateIntp} and for $\delta_3=\frac{pb}{4}$ \eqref{FinalEstimateIntp} into \eqref{SecondStepMainDerivation} so obtaining that for 
\end{proof}
\end{lemma}
With this gained bound, we are in a position to prove our claims.
\section{Proof of the claims}\label{SectionProofsClaims}
\subsection*{Proof of Theorem \ref{MainTheorem}} Bound \eqref{uniformBoundForLpNormSecondCase} ensures the finiteness of the integrals $\id \rho^{p}dx$ on $(0,T_{\max})$ and for any $p>\bar{p}$, through a $p$-dependent constant $K$.  Our goal is to derive the recursive inequality \eqref{RecursiveInequality} of Lemma \ref{RecursiveInequalityLemma} for the sequence
\begin{equation}\label{DefinitionSequelMk}
M_k:=\max\bigg\{1,\sup_{t\in (0,T_{\max})}\id \rho^{p_k}dx\bigg\}\quad \textrm{for all } k\in \N,
\end{equation}
and with $p_k:=2^k+\bar{p}$, introduced in Lemma \ref{Sequel_p_r_theta_Lemma}. 

We start reconsidering inequality \eqref{FirstStepMainDerivation}, which dropping the term $+pb  \id \rho^{p+\alpha-1}dx\id \rho^\beta dx$ reads for $p=p_k$ and for all $t\in (0,T_{\max})$
\begin{equation}\label{FirstStepMainDerivation-Iteration}
\begin{split}
\frac{d}{d t} \id \rho^{p_k} dx&  +2c_1\id \lvert \nabla \rho^\frac{m+p_k-1}{2}\rvert^2dx
\leq p_k a\id \rho^{p_k+\eta-1}dx +n (p_k-1)\alpha_n \id \rho^{p_k+1}dx.
\end{split}
\end{equation}
Successively, and similarly to what done above, we  derive these bounds.
\subsection*{\quad Estimating $n (p_k-1)\alpha_n \id \rho^{p_k+1}$.} Setting $q_k=2\frac{p_k+1}{m+p_k-1}$ and $r_k=\frac{2p_{k-1}}{m+p_k-1}$, we have 
$$1\leq r_k <q_k<l\quad \textrm{and}\quad\frac{q_k}{r_k}<\frac{2}{r_k}+1-\frac{2}{l}.$$ 
Subsequently, by using inequality \eqref{InequalityTipoG-N} with $\gamma_{1,k}$ defined in Lemma \ref{Sequel_p_r_theta_Lemma}, we get on $(0,T_{\max})$ that for every $\tilde{\varepsilon}_1>0$ and a computable positive $\tilde{C}_1(\tilde{\varepsilon}_1)$
\begin{equation}\label{FirstEtimateIntp+1-Iteration}
\begin{split}
n (p_k-1)\alpha_n  \id \rho^{p_k+1}dx&=n (p_k-1)\alpha_n \lVert \rho^\frac{q(m+p_k-1)}{2}\rVert_{L^q(\R^n)}^q\leq \tilde{\varepsilon}_1\id|\nabla \rho^\frac{m+p_k-1}{2}|^2 dx+\tilde{C}_{1,k}(\tilde{\varepsilon}_1)\bigg(\id \rho^{p_{k-1}}dx\bigg)^\frac{\gamma_{1,k}}{r_k}.
\end{split}
\end{equation} 
We herein have to importantly observe that, up to a $k$-independent positive constant $b_1$, it can be shown that
\begin{equation*}%\label{C_Tilde_1_Lambda}
\tilde{C}_{1,k}(\tilde{\varepsilon}_1)=b_1  {{\varepsilon}_1}^\frac{-\lambda_{1,k}q_k}{2-\lambda_{1,k} q_k},
\end{equation*}
where $\lambda_{1,k}$ is as well mentioned in Lemma \ref{Sequel_p_r_theta_Lemma}. In particular, since $\lambda_{1,k} \nearrow 0$
\begin{equation}\label{LimitOfC_Tilde_1_Lambda}
\lim_{k\rightarrow\infty}\tilde{C}_{1,k}(\tilde{\varepsilon}_1)=b_1.
\end{equation}
\subsection*{Estimating $p_k a\id \rho^{p_k+\eta-1}dx$.} For $q_k=2\frac{\eta +p_k-1}{m+p_k-1}$ and  $r_k=\frac{2p_{k-1}}{m+p_k-1}$ we get 
$$1\leq r_k <q_k<l\quad \textrm{and}\quad\frac{q_k}{r_k}<\frac{2}{r_k}+1-\frac{2}{l},$$ so having from \eqref{InequalityTipoG-N} that for all  $\tilde{\varepsilon}_2>0$, some positive $\tilde{C}_{2,k}(\tilde{\varepsilon}_2)$ and $\gamma_{\eta,k}$ as in Lemma \ref{Sequel_p_r_theta_Lemma} 
\begin{equation}\label{FirstEtimateIntp+eta+1-Iteration}
\begin{split}
p_k a \id \rho^{\eta +p_k-1}&=p_k a  \lVert \rho^\frac{q(m+p_k-1)}{2}\rVert_{L^q(\R^n)}^q\leq \tilde{\varepsilon}_2\id|\nabla \rho^\frac{m+p_k-1}{2}|^2 dx +\tilde{C}_{2,k}(\tilde{\varepsilon}_2)\bigg(\id \rho^{p_{k-1}}dx\bigg)^\frac{\gamma_{\eta,k}}{r_k}\quad \textrm{for all } t\in (0,T_{\max}).
\end{split}
\end{equation} 
Reasoning as in the previous case, for some $k$-independent constant $b_2>0$ we have
\begin{equation}\label{LimitOfC_Tilde_2_Lambda}
\lim_{k\rightarrow\infty}\tilde{C}_{2,k}(\tilde{\varepsilon}_1)=b_2\quad \textrm{where }\; \tilde{C}_{2,k}(\tilde{\varepsilon}_2)=b_2  {{\varepsilon}_2}^\frac{-\lambda_{\eta,k}q_k}{2-\lambda_{\eta,k} q_k}.
\end{equation}
\subsection*{Estimating $\id \rho^{p_k}dx$.}
We now take $q_k=\frac{2p_k}{m+p_k-1}$ and  $r_k=\frac{2p_{k-1}}{m+p_k-1}$, and we obtain
$$1\leq r_k <q_k<l\quad \textrm{and}\quad\frac{q_k}{r_k}<\frac{2}{r_k}+1-\frac{2}{l}.$$ In this way, using once again the Gagliardo--Nirenberg-type inequality \eqref{InequalityTipoG-N} and  Lemma \ref{Sequel_p_r_theta_Lemma} with $\gamma_{0,k}$ therein defined, for every $\tilde{\varepsilon}_3>0$ and some $b_3>0$ the sequence
\begin{equation}\label{LimitOfC_Tilde_3_Lambda}
\tilde{C}_{3,k}(\tilde{\varepsilon}_3)=b_3  {{\varepsilon}_3}^\frac{-\lambda_{0,k}q_k}{2-\lambda_{0,k} q_k}\quad \textrm{satisfying }\;\lim_{k\rightarrow\infty}\tilde{C}_{3,k}(\tilde{\varepsilon}_3)=b_3,
\end{equation}
has the property that 
\begin{equation}\label{FirstEtimateIntp-Iteration}
\begin{split}
\id \rho^{p_k}dx&=\lVert \rho^\frac{q(m+p_k-1)}{2}\rVert_{L^q(\R^n)}^q\leq \tilde{\varepsilon}_3\id |\nabla \rho^\frac{m+p_k-1}{2}|^2 dx +\tilde{C}_{3,k}(\tilde{\varepsilon}_3)\bigg(\id \rho^{p_{k-1}}dx\bigg)^\frac{\gamma_{0,k}}{r_k}\quad \textrm{on } (0,T_{\max}).
\end{split}
\end{equation} 
Next we take $\tilde{\varepsilon}_1=\tilde{\varepsilon}_2=\frac{c_1}{2}$, 
and $\tilde{\varepsilon}_3=c_1$, introduce  ${\mu_{1,k}}, {\mu_{\eta,k}}$ and ${\mu_{0,k}}$ from Lemma \ref{Sequel_p_r_theta_Lemma},  and plug inequalities \eqref{FirstEtimateIntp+1-Iteration}, \eqref{FirstEtimateIntp+eta+1-Iteration}, \eqref{FirstEtimateIntp-Iteration} into estimate \eqref{FirstStepMainDerivation-Iteration};  we see that 
\begin{equation*}
\begin{split}
\frac{d}{d t} \id \rho^{p_k} dx& \leq -\id \rho^{p_k} dx+\tilde{C}_{1,k}\big(\frac{c_1}{2}\big)\Big(\id\rho^{p_{k-1}}dx\Big)^{\mu_{1,k}}+\tilde{C}_{2,k}\big(\frac{c_1}{2}\big)\Big(\id\rho^{p_{k-1}}dx\Big)^{\mu_{\eta,k}}\\&\quad +\tilde{C}_{3,k}(c_1)\Big(\id\rho^{p_{k-1}}dx\Big)^{\mu_{0,k}}\quad \textrm{ on } (0,T_{\max}).
\end{split}
\end{equation*}
At this juncture, we define 
\begin{equation}\label{HkDefinition}
H_k:=\tilde{C}_{1,k}\big(\frac{c_1}{2}\big)+\tilde{C}_{2,k}\big(\frac{c_1}{2}\big)+\tilde{C}_{3,k}(c_1)\quad \textrm{for all } \,k\in \N,
\end{equation}
and recalling \eqref{DefinitionSequelMk} and taking in mind the bounds of Lemma \ref{Sequel_p_r_theta_Lemma} we obtain this initial value problem 
\begin{equation*}
\begin{cases}
\frac{d}{d t} \id \rho^{p_k} dx\leq -\id \rho^{p_k} dx+H_kM_{k-1}^{2} \quad \textrm{ on } (0,T_{\max}),\\
\id \rho^{p_k}(x,0) dx=\id \rho_0^{p_k} dx,
\end{cases}
\end{equation*}
which implies 
\begin{equation*}%\label{EstimatingMkWithHK} 
M_k\leq \max\Big\{\id \rho_0^{p_k} dx,H_kM_{k-1}^{2}\Big\}  \quad \textrm{ for all } k\in \N.
\end{equation*}
From the convergence of the sequences $\tilde{C}_{i,k}(\cdot)$ for $i=1,2,3$ (indeed \eqref{LimitOfC_Tilde_1_Lambda}, \eqref{LimitOfC_Tilde_2_Lambda} and \eqref{LimitOfC_Tilde_3_Lambda} hold), and the definition of $H_k$ in \eqref{HkDefinition}, we have that for some $L>0$ it is attained from the previous relation that 
\begin{equation}\label{EstimatingMkWithHK}
M_k\leq L M_{k-1}^{2}  \quad \textrm{ for all } k\in \N.
\end{equation}
In this circumstances,  we have to distinguish two cases.
\begin{itemize}
\item If $LM_{k-1}^{2}<\id \rho_0^{p_k}dx$, then
\begin{equation*}
\begin{split}
M_{k-1}^{\frac{1}{p_{k-1}}}&=\sup_{t\in(0,T_{\max})}\bigg(\id \rho^{p_{k-1}}dx\bigg)^\frac{1}{p_{k-1}}<\bigg(\frac{1}{L}\id \rho_0^{p_{k}}dx\bigg)^\frac{1}{2p_{k-1}}
=\frac{1}{(L)^\frac{1}{2 p_{k-1}}}\Big[\bigg(\id \rho_0^{p_{k}}dx\bigg)^\frac{1}{p_{k}}\Big]^\frac{p_k}{2 p_{k-1}},
\end{split}
\end{equation*}
so that since we have $\frac{1}{2 p_{k-1}} \nearrow 0$ and $\frac{p_k}{2p_{k-1}}  \nearrow 1$, we can conclude that
\begin{equation}\label{BoundIterative1}
\begin{split}
\lVert \rho (\cdot, t)\rVert_{L^\infty(\R^n)}&:=\lim_{k\rightarrow \infty}\sup_{t\in (0,T_{\max})}\bigg(\id \rho^{p_{k-1}}dx\bigg)^\frac{1}{p_{k-1}} =\lim_{k\rightarrow \infty} M_{k-1}^{\frac{1}{p_{k-1}}}< \lVert \rho_0 (\cdot)\rVert_{L^\infty(\R^n)} \quad \textrm{for all } t\in  (0,T_{\max}).
\end{split}
\end{equation}
\item If $LM_{k-1}^{2}\geq \id \rho_0^{p_k}dx$, then
 from inequality \eqref{EstimatingMkWithHK} an application of Lemma \ref{RecursiveInequalityLemma} infers
\begin{equation*}
M_k\leq L^{k}\cdot M_{0}^{ 2^k}  \quad \textrm{ for all } k\in \N,
\end{equation*}
and since $\frac{k}{p_k} \nearrow 0$ and $\frac{2^{k}}{p_k} \nearrow 1$, it holds that 
\begin{equation}\label{BoundIterative2}
\begin{split}
\lVert \rho (\cdot, t)\rVert_{L^\infty(\R^n)}&:=\lim_{k\rightarrow \infty}\sup_{t\in (0,T_{\max})}\bigg(\id \rho^{p_{k}}dx\bigg)^\frac{1}{p_{k}} =\lim_{k\rightarrow \infty} M_{k}^{\frac{1}{p_{k}}}\leq \lim_{k\rightarrow \infty} L^\frac{k}{p_k}\cdot M_{0}^{{\frac{2^{k}}{p_k}}}=M_0\quad \textrm{on }  (0,T_{\max}).
\end{split}
\end{equation}
From the last estimates \eqref{BoundIterative1} and \eqref{BoundIterative2} we get
$$ 
\lVert \rho (\cdot, t)\rVert_{L^\infty(\R^n)}\leq C:=\max\{\lVert \rho_0 (\cdot)\rVert_{L^\infty(\R^n)},M_0\},
$$
and the extensibility criterion \eqref{ExtensibilityCrit} provides $T_{\max}=\infty.$
\qed
\end{itemize}
\subsection*{Proof of Corollary \ref{MainTheoremCorollary}}We have only to observe that in the repulsive case the sign ``$+$'' in estimate \eqref{FirstStepMainDerivation} associated to term $+n (p-1)\alpha_n \id \rho^{p+1}dx$  is ``$-$'', so we can neglect it. Henceforth, retracing the proof of Lemma \ref{LEmmaBoundinlP}, we see that the introduction of $\tilde{\Lambda}_1$, which leads to the condition $\alpha+\beta>\frac{2\eta l -ml-2\eta}{l-2}$, is superfluous. Consequently, assumptions \eqref{b1} suffice to conclude similarly to what already explained.
\qed
\subsubsection*{Acknowledgments}
%The author is grateful to the referees for helpful suggestions which improved this article.
GV is a member of the Gruppo Nazionale per l'Analisi Matematica, la Probabilit\`a e le loro Applicazioni (GNAMPA) of the Istituto Na\-zio\-na\-le di Alta Matematica (INdAM) and  is  partially supported by the research projects \textit{Integro-differential Equations and Non-Local Problems}, funded by Fondazione di Sardegna (2017), and by MIUR (Italian Ministry of Education, University and Research) Prin 2017 \textit{Nonlinear Differential Problems via Variational, Topological and Set-valued Methods} (Grant Number: 2017AYM8XW). The research of TL is supported by NNSF of P. R. China (Grant No. 61503171), CPSF (Grant No. 2015M582091), and NSF of Shandong Province (Grant No. ZR2016JL021).

% % % % % % %
%\bibliography{Bibliography}{}
%\bibliography{reference}
\bibliographystyle{abbrv}

\end{document}